\documentclass{amsart}
\usepackage[english]{babel}
\usepackage{latexsym}
\usepackage{amsmath,amssymb,amsfonts,amsthm}
\usepackage{mathrsfs}
\usepackage{epsfig}
\usepackage[all]{xy}
\usepackage{hyperref}

\title[Determinantal Classification of Foliations]{Determinantal Classification of Foliations on $\CP$ with a Unique Singular Point with Nonzero Linear Part}
\author{Claudia R. Alc\'antara}
\date{}
\address{Departamento de Matem\'aticas, Universidad de Guanajuato.
 Callej\'on Jalisco s/n, A.P. 402, C.P. 36000,
Guanajuato, Gto. M\'exico.}
\email{claudia@cimat.mx, ce.reynoso@ugto.mx}
\subjclass[2010]{Primary 32S65, 14M12. Secondary: 37F75.}
 \keywords{holomorphic foliation, singularity, singular scheme, determinantal variety.}
 
 \newtheorem{theorem}{Theorem}
\newtheorem{lem}{Lemma}
\newtheorem{example}{Example}
\newtheorem{defin}{Definition}
\newtheorem{prop}{Proposition}

\newtheorem{algorithm}{Algorithm}
\newtheorem{remark}{Remark}
\newcommand{\Oc}{\mathcal{O}}

\newcommand{\CP}{\mathbb{CP}^2}
\newcommand{\C}{\mathbb{C}}

\newcommand{\N}{\mathbb{N}}
\newcommand{\F}{\mathcal{F}}

\begin{document}
\maketitle

\begin{abstract}
We prove that the space of holomorphic foliations on $\CP$ of degree $d \geq 3$, with a unique singular point with non-zero linear part
is a determinantal variety. We provide an algorithm for describing the equations that define it. We give explicit equations for $d=3$, and use this result to construct an explicit family that generalizes a previously known one. For degree $4$, we exhibit an example of a foliation satisfying these conditions (no examples of this type were previously known). Finally, we prove that the singular scheme of any foliation of degree $d$ with a unique singular point with non-zero linear part is contained in a smooth curve of order strictly less than $d$.
\end{abstract}

\tableofcontents

\section{Introduction}

Although holomorphic foliations on the complex projective plane have been extensively studied, many classification problems remain open. Among them is the problem of classifying foliations on $\CP$ of degree $d$ with a unique singular point, that is, foliations having a single singularity with Milnor number $d^2+d+1$. The study of this class of foliations has attracted increasing attention in recent years. Some papers that study foliations on $\CP$ with a unique singularity are:  \cite{alcantara2}, \cite{alcantara3}, \cite{alcantara-mozo}, \cite{alcantara-zamora}, \cite{castorena}, \cite{cerveau}, \cite{cerveaudeserti}, \cite{coutinho}, \cite{coutinho1}, \cite{coutinhojales}, \cite{coutinhooliveira}, \cite{coutinhooliveira2} and \cite{fernandez}.

One of the main motivations for studying this class of foliations is that it provides a natural framework for investigating several fundamental problems in the theory of holomorphic foliations. The reason is that the singular scheme is supported at a single point, while its scheme-theoretic structure may still be highly nontrivial. This allows us to use powerful tools from the theory of zero-dimensional schemes, which are particularly well suited to studying schemes supported at a single point. Some of the fundamental problems we have in mind are the following.

A first problem concerns the construction of foliations without algebraic leaves. Although it is known that foliations without algebraic invariant curves form a dense nonempty subset of the space of foliations, explicit examples are remarkably difficult to construct (see for example \cite{{joanoulou}} and \cite{zoladek} for foliation on the projective plane and \cite{lizarbe} for foliations on weighted projective planes). In this direction, it was proved in \cite{alcantaraPantaleon} that every foliation on $\CP$ with a unique singularity having non-zero linear part has no algebraic leaves. Consequently, a classification of foliations with a unique singularity would provide a systematic source of explicit examples of foliations without algebraic invariant curves.

A second motivation is related to the classical Poincaré problem. Given a foliation on $\CP$ admitting a rational first integral, one asks whether the degree of such an integral can be bounded just in terms of the degree of the foliation. It is now known that the answer is negative, although only a few families of counterexamples are currently available. Recently, Alcántara and Zamora  in \cite{alcantara-zamora} constructed a family of foliations with a unique singularity whose rational first integrals have unbounded degree. This suggests that the geometry of the singular scheme, which in this case is supported at a single point, may provide new insight into the obstructions preventing the existence of such a bound.

A third important motivation comes from the study of the Exceptional Minimal Set Conjecture, which states that for a foliation on $\CP$ there are no exceptional minimal sets disjoint from the singular set (see \cite{camacho}). In this context, the classification of foliations on $\CP$ with a unique singularity may provide a suitable setting for understanding this problem, since the highly constrained geometry of the singular scheme could reveal new phenomena that are not apparent in the general case. 

A fourth motivation comes from the problem of determining a foliation from its singular scheme. This question has been extensively studied by Campillo and Olivares, who proved that, under suitable assumptions, there exist proper subschemes of the singular scheme satisfying certain geometric properties that still determine the foliation (see \cite{campillo2}). In their work on Cayley--Bacharach properties of singular schemes (see \cite{campillo3}), they exhibited a foliation with a unique singularity for which no such subscheme exists. This shows that foliations with a unique singularity already contain subtle phenomena related to the geometry of singular schemes, providing further motivation for their systematic study.

Therefore, the main goal of this work is to provide a classification, that is as explicit and constructive as possible, of foliations on $\CP$ of degree greater than or equal to 3  that have a unique singularity having non-zero linear part, meaning that the singularity is either of saddle-node type or nilpotent type.
\\

The complete classification of foliations with a unique singular point is known only for degrees $1$ and $2$. In \cite{cerveau} we can find the proof that, up to change of coordinates, there are only four foliations
with a unique singularity, they are: 

\begin{align*}
\omega_1=&-z^3\,dx+y^2z\,dy
+(xz^2-y^3)\,dz,\\
\omega_2=&z^3\,dx-(y^2z+xz^2)\,dy
+(xyz-xz^2+y^3)\,dz,\\
\omega_3=&-yz^2\,dx+(y^2z+z^3)\,dy
+(xyz-y^3-yz^2)\,dz,\\
\omega_4=&(y^3-xyz-z^3)\,dx
-(xy^2+yz^2)\,dy+(x^2y+xz^2+y^2z)\,dz.
\end{align*}

The multiplicity of the unique singular point in $\omega_1, \omega_2$ and $\omega_3$ is $2$, while in the last one we have a unique singularity of type saddle-node. Therefore, there are no nilpotent singularities with Milnor number $7$ for foliations of degree $2$.
\\

In \cite{alcantara-ronzon} we have some families of foliations of degree $3$ with a unique singularity (Milnor number $13$), and with multiplicity $2$ and $3$. In \cite{alcantara1} we describe a locally closed algebraic variety parametrizing foliations on $\CP$ of degree $d$ with a unique singularity with multiplicity $d$. The techniques used in these works come from representation theory. We use the irreducible representation given by change of coordinates in the space of foliations on $\CP$ of degree $d$. With this tools we can obtain foliations with a unique singular point with multiplicity greater that $\frac{d}{2}$. For multiplicity one, or equivalently, when the linear part at the singularity is nonzero, the algebraic techniques are based on Groebner bases. We briefly explain this approach below.
\\

Using Gr\"{o}bner basis with respect to the lexicographical monomial order we obtain necessary and sufficient conditions to know, given a foliation on $\CP$ of degree $d$, if it has a singular point with Milnor number $d^2+d+1$ and with non-zero linear part (see Theorem 1 of \cite{alcantara2}). The singular schemes of these foliations are curvilinear schemes  (i.e. the scheme is on a smooth curve ) with support in one point $p$. The set of curvilinear schemes $\mathbb{H}_{curv,p}^{d^2+d+1}(\CP)$ is a subvariety of the Hilbert scheme of $d^2+d+1$ points in $\CP$, which we denote by $\mathbb{H}^{d^2+d+1}(\CP)$. The variety $\mathbb{H}_{curv,p}^{d^2+d+1}(\CP)$  is irreducible, non-singular of dimension $d^2+d$ (see  \cite{gottsche}). 
Therefore, the goal of this work  is to describe the subvariety  of $\mathbb{H}_{curv,p}^{d^2+d+1}(\CP)$ that parametrizes foliations on the projective plane of degree $d$ with a unique singular point with non-zero linear part.
\\

Thus, relying on properties of the curvilinear scheme and on elementary linear algebra, we obtain the following results. Our first main result shows that the space of foliations on $\CP$ of degree $d$ with a unique singular point with non-zero linear part  has a determinantal structure (see proposition \ref{determinantal}). More precisely, we provide a general procedure to construct the matrix whose rank conditions yield the determinantal equations defining this space. This description is valid for every degree \(d\geq 3\) and provides an effective approach to studying these foliations through explicit algebraic equations.

Next, we investigate the local structure of the singular scheme of these foliations. We prove that, in suitable local coordinates, this scheme is defined by an ideal of the form
\[
\left \langle y-\sum_{i=2}^{d^2+d}a_i z^i,\ z^{d^2+d+1}\right \rangle \subset \C[y,z],
\]
where $\min\{i\mid a_i\neq 0\}<d$ (see Theorem \ref{multiplicity}). In particular, this improves the previously known bound (\(\min\{i\mid a_i\neq 0\}\leq d\) proven in Theorem 2 of \cite{alcantara2}), by showing that the inequality is always strict. Thus, the local singular scheme of a foliation  on $\CP$ of degree $d$ is supported on a smooth curve whose order is strictly less than \(d\).

We then apply the above results to degrees $3$ and $4$. In degree \(3\), we obtain explicit determinantal equations describing the corresponding space (see Theorem \ref{degree3}). These equations allows us to construct an explicit family of foliations that generalizes a previously known family (see example \ref{familydegree3}). We provide the first example of a foliation on $\CP$ of degree $4$ with a unique singular point having nonzero linear part. In this case, the singularity is a saddle-node (see proposition \ref{exampledegree4}). This is an example that several authors had been looking for. 

It is worth emphasizing that all the foliations found in these varieties are foliations without algebraic leaves, the proof of this is Theorem 3 of \cite{alcantara3}. As we know, although a generic foliation on $\CP$ has no algebraic leaves, finding explicit examples of such foliations is a very difficult task. Thus, these methods provide us with a large number of such examples.
\\

The article is organized as follows. Section 2 collects the preliminaries on foliations on $\CP$, introduces the notion of singular scheme of a foliation, and recalls some basic facts about curvilinear schemes and Gröbner bases that will be needed throughout the paper. Section 3 is devoted to the main results of this work.

\section{Preliminaries}

\subsection{Holomorphic foliations on $\CP$}

In this section we give the basic definitions and results that we need about holomorphic foliations on $\CP$ of degree $d$ to develop the work. 
\\

 Up to multiplication by a nonzero scalar there are two equivalent ways to describe a \textbf{foliation on $\CP$
of degree $d$} (see \cite{gmont}):
\\

\begin{enumerate}

\item By a homogeneous vector field:
\begin{equation*}
X=P(x,y,z) \frac{\partial}{\partial x}+ Q(x,y,z) \frac{\partial}{\partial y}+ R(x,y,z) \frac{\partial}{\partial z}=
{\left(\begin{array}{ccc}
P(x,y,z) \\
Q(x,y,z)\\
R(x,y,z)
\end{array}\right)}
\end{equation*}

\noindent where $P,Q,R \in \C[x,y,z]$ are homogeneous of degree $d$. And if we consider the radial foliation

\begin{equation*}
E=x \frac{\partial}{\partial x}+ y \frac{\partial}{\partial y}+ z \frac{\partial}{\partial z},
\end{equation*}

\noindent then $X$ and $X+F(x,y,z)E$ represent the same foliation for all $F \in \C[x,y,z]$ homogeneous of degree $d-1$. 
\\

\item By a homogeneous 1-form:  $\omega= A(x,y,z) dx+ B(x,y,z)dy+C(x,y,z)dz$, such that $A,B,C \in \C[x,y,z]$ are homogeneous of degree $d+1$ and these satisfy 
the Euler's condition $xA+yB+zC=0$.
\end{enumerate}

If we have a foliation on $\CP$
of degree $d$ given by the homogeneous vector field $P(x,y,z) \frac{\partial}{\partial x}+ Q(x,y,z) \frac{\partial}{\partial y}+ R(x,y,z) \frac{\partial}{\partial z}$, then the corresponding 1-form is $(zQ-yR)dx-(zP-xR)dy+(yP-xQ)dz$.

\vspace{.3cm}

With this we can see that the space $\mathcal{F}_d$ of foliations on $\CP$ of degree $d$ is a projective space of dimension $d^2+4d+2$. 
Now we define the notion of singular point for a foliation and two important invariants for this.

\begin{defin} The point $p=(a:b:c) \in \CP$ is \textbf{singular} for the foliation $\omega= A(x,y,z) dx+ B(x,y,z)dy+C(x,y,z)dz$ if 
  
\begin{equation*}
A(a,b,c)=B(a,b,c)=C(a,b,c)=0.
\end{equation*}

\noindent The set of singular points of $\omega$ will be denoted by $Sing(\omega)$.
\end{defin}

\begin{defin} Let $f(y,z)dy+g(y,z)dz$ be a local generator of $\omega$ in $p=(1:b:c)$. Then
\vspace{.2cm}

\noindent \textbf{the Milnor number of $p$} is $\mu_p(\omega):=\dim_{\C}\frac{\Oc_{\C^2,(b,c)}}{<f,g>}$, \\
\noindent \textbf{the multiplicity of $p$} is $m_p(\omega):=min\{ord_{(b,c)}(f),ord_{(b,c)}(g)\}$.
\\

The multiplicity of $p$ is $1$ if and only if we have non-zero linear part. Suppose that $m_p(\omega)=1$, then if we have one proper value zero and one proper value non-zero the singularity is called \textbf{saddle-node}. And it is \textbf{nilpotent} if both proper values are zero.
\end{defin}

As we can see in the following result, the number of singularities depends only on the degree of the foliation. 

\begin{prop} (see \cite{brunella}) Let $X$ be a foliation on $\CP$ of degree $d$ with isolated singularities then

\begin{equation*}
d^2+d+1=\sum_{p \in \CP} \mu_p(X).
\end{equation*}
\end{prop}

\subsection{The Singular Scheme of a Foliation on $\CP$ via Gröbner Bases}

The well-known Hilbert scheme $\mathbb{H}^{n}(\CP)$, of $n$ points on $\CP$ is a space which parametrizes zero-dimensional subschemes of $\CP$ of  length $n$, or 
homogeneous saturated ideals $J \in \C[x,y,z]$ such that the dimension as $\C$-vector space of $\frac{\C[x,y,z]}{J}$ is $n$. 
This is a smooth, irreducible algebraic variety of dimension $2n$ (see \cite{fogarty}). 
\\

In a coordinate chart of $\CP$ we can see a zero-dimensional subscheme of  length $n$
 as $n$ points in the plane $\C^2$, the space of these subschemes is denoted by $\mathbb{H}^{n}(\C^2)$. The Hilbert scheme $\mathbb{H}^{n}(\C^2)$  can be identified with the set of ideals $I$ in $\C[y,z]$, such that the dimension as $\C$-vector space of $\frac{\C[y,z]}{I}$ is $n$. \\

Now we will define the singular scheme of a foliation on $\CP$ of degree $d$ with isolated singularities given by the 1-form, $\omega= A(x,y,z) dx+ B(x,y,z)dy+C(x,y,z)dz$.

\begin{defin}  The \textbf{singular scheme $\mathcal{S}ing(\omega)$} of a foliation on $\CP$ of degree $d$ with isolated 
singularities, defined by the 1-form $\omega=A(x,y,z) dx+ B(x,y,z)dy+C(x,y,z)dz$, is given by the homogeneous ideal:

\begin{equation*}
I=\Big \langle A(x,y,z), B(x,y,z), C(x,y,z) \Big\rangle \subset \C[x,y,z].
\end{equation*} 

\noindent This ideal satisfies $\dim_{\C} \frac{\C[x,y,z]}{I}=d^2+d+1$, therefore is a point in $\mathbb{H}^{d^2+d+1}(\CP)$.
 \end{defin}

As we have just seen, the singular scheme of a foliation on $\CP$ of degree $d$ with isolated singularities is a point in the Hilbert scheme 
$\mathbb{H}^{d^2+d+1}(\CP)$.
Moreover, a foliation of degree greater or equal than $2$ is completely determined by its singular scheme. We can see the proof of this in Theorem 3.5 of \cite{campillo1}.
Therefore, we can approach the classification of foliations  on $\CP$ of degree $d\geq 2$ by studying their singular schemes. To carry out the necessary computations, we will work in an affine coordinate chart.
\\

Let $\omega=A(x,y,z) dx+ B(x,y,z)dy+C(x,y,z)dz$ be a foliation on $\CP$ of degree $d$ with isolated singularities. With a change of coordinates we can move all the singularities to the chart $U_0=\{(1:y:z): (y,z) \in \C^2\}$, the foliation in this chart is:

\begin{align*}
&B(1,y,z)dy+C(1,y,z)dz\\  
\end{align*}

\noindent therefore the local singular scheme is given by the ideal $\langle B(1,y,z), C(1,y,z) \rangle$ in $\C[y,z]$, which is a point in the Hilbert scheme $\mathbb{H}^{d^2+d+1}(\C^2)$. 
\\

We can study this ideal in two variables using Gr\"{o}bner basis. In this work we  give only the definitions and results that we will need. 
A very good reference to see details on  Gr\"{o}bner basis theory is the book \cite{cox}. 
\\

The first thing we must do when we study Gr\"{o}bner basis is to fix a monomial order, which is a total, multiplicative and well order in the set of monic monomials of a given polynomial ring. In this work we use the following:

\begin{defin} The \textbf{lexicographical monomial order $>_{lex}$ in $\C[y,z]$} is defined as follows:

\begin{equation*}
y^{\alpha_1}z^{\alpha_2} >_{lex} y^{\beta_1}z^{\beta_2}
\end{equation*}

\noindent if and only if if the leftmost non-zero entry of $(\alpha_1,\alpha_2)-(\beta_1,\beta_2)=(\alpha_1-\beta_1,\alpha_2-\beta_2)$ is positive.

\end{defin}

\begin{remark} For this order we have: $y >_{lex} z^n$ for all $n \in \N$ and $z^{\beta_1} >_{lex} z^{\beta_2}$ if and only if $\beta_1 > \beta_2$.
\end{remark}

Now we are going to introduce other important definitions.

\begin{defin}
\begin{enumerate}
\item Let $f(y,z) \in \C[y,z]-\{0\}$, then the \textbf{leading monomial of $f$} is the greater non-zero monomial with respect to $>_{lex}$. It is denoted by $Lm(f)$.

\item Let $I \subset \C[y,z]$ be a non-zero ideal, we define the \textbf{initial ideal of $I$ with respect to $>_{lex}$} as $in(I)=\langle Lm(f): f \in I \rangle$.
\end{enumerate}
\end{defin}

Finally we will give the definition of Gr\"{o}bner basis for an ideal.

\begin{defin} Let $I$ be a non-zero ideal in $\C[y,z]$. The set $\mathcal{G}=\{f_1,...,f_s\} \subset \C[y,z]$ is a Gr\"{o}bner basis for $I$ if and only if 
$in(I)=\langle Lm(f_i): f_i \in \mathcal{G} \rangle$. 
\\

\noindent A Gr\"{o}bner basis $\mathcal{G}$ for $I$ is called \textbf{reduced} if and only if for all
$j=1,...,s$ the polynomial $f_j$ is monic and its leading monomial $Lm(f_j)$ does not divide any term of $f_i$ for all $i \neq j$.

\end{defin}

Something easy to prove and very important in this theory is the following (corollary 2.5.6 in \cite{cox}).

\begin{prop} Let $I$ be a non-zero ideal in $\C[y,z]$. If $\mathcal{G}$ is a Gr\"{o}bner basis for $I$ then $I$ is generated by $\mathcal{G}$.
\end{prop}

There exists an algorithm due to Buchberger for computing a Gr\"{o}bner basis for any non-zero ideal in a polynomial ring. The following proposition talks about the uniqueness of the reduced Gr\"{o}bner basis. The proof can be found in
proposition 2.7.6 of \cite{cox}.

\begin{prop} Let $I \subset \C[y,z]$ be an ideal such that $I \neq \{ 0\}$. Then $I$ has a unique reduced Gr\"{o}bner basis with respect to the lexicographical monomial order.
\end{prop}

Of course these definitions are valid for any monomial order, but in this paper we only work with $>_{lex}$. 
\\

Let $p \in \C^2$, we will denote $\mathbb{H}_{p}^{n}(\C^2)=\big\{I \in \mathbb{H}^{n}(\C^2): \mathbb{V}(I)=\{p\} \big\}$, as the set of zero-dimensional subschemes of $\C^2$ of  length $n$, with support at $p$. This is an irreducible, non-singular algebraic variety of dimension $n-1$ (see  \cite{briancon}). A curvilinear scheme of  length $n$ with support at $p$  is a scheme on a smooth curve. Finally we have that the set of curvilinear schemes $\mathbb{H}_{curv,p}^{n}(\C^2)$ with support at $p$, is an open  irreducible, non-singular subvariety of the Hilbert scheme $\mathbb{H}_p^{n}(\C^2)$ (see  \cite{gottsche}). If we fix the point $p=(0,0)$ then, up to change of coordinates, and ideal $I \in \mathbb{H}_{curv,p}^{n}(\C^2)$ has the form $\langle y+\sum_{i=2}^{n-1}a_i z^i, z^{n}\rangle$.

\section{Determinantal variety of Foliations on $\CP$ of degree $d$ with a unique singular point}

In this section, we develop the construction of the determinantal variety that parametrizes foliations on $\CP$ of degree $d$ with a unique singularity with multiplicity $1$ (or non-zero linear part). For thar we work at the chart $U_0=\{(1:y:z) \in \CP\}$, we fix the singularity at $(0,0)$ and the linear part (up to change of coordinates) as $ydy$ when we have a sadle-node singularity, or as $ydz$ when the singularity is nilpotent. 
\\

Let $\omega$ be a foliation on $\CP$ of degree $d$ and let $f(y,z)dy+g(y,z)dz$ be the affine form at the chart $U_0$.  Theorem 1 of \cite{alcantara2} says that the point $(0,0)$ is a singular point for $\omega$ with Milnor number $d^2+d+1$ and the linear part around the point is $ydy$ (when the singularity is a sadle-node) or $ydz$ (when the singularity is nilpotent), if and only if, the reduced Gr\"{o}bner basis of the ideal $\big \langle f(y,z), g(y,z) \big \rangle \subset \C[y,z]$, with respect to the lexicographical monomial order where $y > z$, has the curvilinear form:

$$I=\Big\langle y-\sum_{i=2}^{d^2+d} a_iz^i, z^{d^2+d+1} \Big\rangle \subset \mathbb{C}[y,z],$$

\noindent for some $(a_2,...,a_{d^2+d}) \in \mathbb{A}^{d^2+d-1}$. Therefore, the main objective of this work is to study the algebraic variety

\begin{align*}
\F_d^0:=\Big\{&(a_2,...,a_{d^2+d}) \in \mathbb{A}^{d^2+d-1} : \Big\langle y-\sum_{i=2}^{d^2+d} a_iz^i, z^{d^2+d+1} \Big\rangle\\
  &\textrm{ is the local singular scheme of a foliation on $\CP$ of degree $d$} \Big\},
\end{align*}

\noindent That is, we aim to find polynomial equations defining this subvariety of the affine space $\mathbb{A}^{d^2+d-1} $. We will show that these equations can be obtained explicitly, using results mainly from linear algebra. We begin by stating some useful lemmas that will help us establish our main results.

\begin{lem} \label{lemma1}  The ideal
\[
I=\left\langle y-\sum_{i=2}^{d^2+d}a_i z^i,\,
z^{d^2+d+1}\right\rangle \subset \C[y,z]
\]
defines the local singular scheme of a foliation on \(\mathbb{CP}^2\) of degree $d$ in the affine chart $U_0$, if and only if there exist three polynomials \(A,B,C\in\mathbb{C}[y,z]\), of degree \(d+1\), such that \(A,B,C\in I\), and the ideal
$\left\langle y-\sum_{i=2}^{d^2+d-1}a_i z^i,\,
z^{(d-j)(d+j+1)}\right\rangle$ contains no nonzero polynomial of degree \(d-j\), for $0 \leq j < d$.\end{lem}

\begin{proof}
The proof follows directly from Theorem 4.7 of \cite{campillo2} and the fact that the scheme defined by  $I$ contains a unique subscheme of length $(d-j)(d+j+1)$ for $0 \leq j < d$, which is defined by the ideal:
$$\left\langle y-\sum_{i=2}^{d^2+d-1}a_i z^i,
z^{(d-j)(d+j+1)}\right\rangle,$$

\noindent see \cite{gottsche}.
\end{proof}

\begin{lem} \label{lemma2}

Let $g(z)\in \mathbb{C}[z]$ and $n \in \N$. Consider the ideal

$$I =  \langle y-g(z), z^n \rangle \subset \mathbb{C}[y,z].$$

\noindent Then for $f(y,z)\in \mathbb{C}[y,z]$, one has that: $f(y,z) \in I$ if and only if, $f(g(z),z)\in \langle z^n \rangle \subset \mathbb{C}[z].$
\end{lem}

\begin{proof} Suppose that $f(y,z) \in I$. Then there exist polynomials $A(y,z), B(y,z)\in \mathbb{C}[y,z]$ such that $f(y,z)=(y+g(z))A(y,z)+z^n B(y,z)$.
Evaluating at $y=g(z)$, the first term vanishes and we obtain $f(g(z),z)=z^n B(g(z),z)$. Hence $f(g(z),z)$ is a multiple of $z^n$.
\\
Conversely, assume that $f(g(z),z)\in  \langle z^n \rangle$, so there exists $h(z)\in \mathbb{C}[z]$ such that $f(g(z),z)=z^n h(z)$. If we consider $f(y,z)\in \mathbb{C}[z][y]$, by the division algorithm, we obtain:

$$f(y,z)=(y+g(z))Q(y,z)+R(z),$$

\noindent where $Q(y,z)\in \mathbb{C}[y,z]$ and $R(z)\in \mathbb{C}[z]$. Evaluating at $y=g(z)$, we get $R(z)=f(g(z),z) \in  \langle z^n \rangle$. So $R(z)=z^n h(z)$ and substituting back, we obtain

$$f(y,z)=(y+g(z))Q(y,z)+z^n h(z),$$

\noindent which shows that $f\in \langle y+g(z), z^n \rangle$.

\end{proof}

 Let $k,n \in \N$, we have that $\mathbb{C}[y,z]_{\leq k}$, the complex vector space of polynomials of degree at most $k$, has dimension $\frac{(k+1)(k+2)}{2}$ and 
 $\frac{\mathbb{C}[z]}{\langle z^{n} \rangle}$ has dimension $n$, over $\C$. 
 
\begin{defin} Given $k,n \in \N$ and $N \in \N_{\geq 2}$, consider
 $p=(a_2,...,a_{N}) \in \mathbb{A}^{N-1}$ we define the linear transformation:

\begin{align*}
\boldsymbol{T(k,n)_p}: \mathbb{C}[y,z]_{\leq k} &\to \frac{\mathbb{C}[z]}{\langle z^{n} \rangle}\\
f(y,z) &\mapsto f\Big(\sum_{i=2}^{N} a_iz^i,z\Big)+\langle z^{n} \rangle.\\
\end{align*}

\end{defin}

A matrix representing $T(k,n)_p$ has size
$n \times \frac{(k+1)(k+2)}{2}$. We denote by \(\boldsymbol{\operatorname{rank}(T(k,n)_p)}\) and \(\boldsymbol{\operatorname{null}(T(k,n)_p)}\) the rank and nullity of the linear transformation \(T(k,n)_p\), respectively.
\\

With this definition in mind, and using lemma \ref{lemma1} and lemma \ref{lemma2}, we obtain the following result. 

\begin{prop} \label{determinantal} Let $d \in \N_{\geq 2}$, then 
\begin{align*}
\F_d^0=\Big\{&p=(a_2,...,a_{d^2+d}) \in \mathbb{A}^{d^2+d-1} : \operatorname{null}( T(d+1,d^2+d+1)_p) \geq 3,\\
& \operatorname{null} (T \big(d-j,(d-j)(d+j+1)\big)_{p})=0 \quad \forall \quad 0 \leq j < d \Big\},
\end{align*}
\end{prop}

\begin{proof} By lemma \ref{lemma2},  a polynomial $f(y,z) \in \mathbb{C}[y,z]_{\leq k}$ is in the ideal $\left\langle y-\sum_{i=2}^{n-1}a_i z^i,
z^{n}\right\rangle$, if and only if $f(y,z) \in Ker\big(T(k,n)_{(a_2,...,a_{n-1})}\big)$. Thus, the conditions stated in lemma \ref{lemma1} translate directly into those given in the statement of this lemma.
\end{proof}

Since $\operatorname{null}( T(d+1,d^2+d+1)_p) \geq 3$ if and only if $rank(T(d+1,d^2+d+1)_p) \leq \frac{d^2+5d}{2}$, we have that  $\F_d^0$ is an open set in the determinantal variety:

$$\mathcal{D}_d= \Big \{p \in \mathbb{A}^{d^2+d-1}: rank\big(T(d+1,d^2+d+1)_p \big) \leq \frac{d^2+5d}{2} \Big\}.$$

By Theorem 3 of \cite{alcantara2}, we know that \(\mathcal{F}_d^0\) is nonempty for every odd \(d>1\). For even degrees, non-emptiness is currently established only for $d=2$, through the example given by Cerveau et al. (see \cite{cerveau}), and for $d=4$, for which we provide an explicit example in the final subsection of this paper.

\begin{algorithm}\label{algorithm} We next describe an algorithm that produces a minimal set of polynomial equations defining $\F_d^0$ and explicitly computes these equations.
\\

 We know that the determinantal variety $\mathcal{D}_d$ is defined by all minors of size $\frac{d^2+5d}{2}+1$ of the $(d^2+d+1)\times\frac{(d+2)(d+3)}{2}$ matrix, which are extremely numerous. But we can describe a much smaller set of equations that will still define the algebraic variety in which $\mathcal{F}_d^0$ is an open subset. If

$$I=\Big\langle y-\sum_{i=2}^{d^2+d}a_iz^i, z^{d^2+d+1} \Big\rangle \subset \C[y,z].$$ 

\noindent is the local singular scheme of a foliation on $\CP$ of degree $d$, then by Theorem  3.2 of \cite{campillo1}, we have that $null(T(d+1,d^2+d+1)_{p=(a_2,...,a_{d^2+d})}=3$. Therefore 
$$rank(T(d+1,d^2+d+1)_{p=(a_2,...,a_{d^2+d})}=\frac{d^2+5d}{2}.$$ 

\noindent We conclude that there exists $A_p$, a $\frac{d^2+5d}{2}\times \frac{d^2+5d}{2}$ block of the matrix $M$ (representing the linear transformation $T(d+1,d^2+d+1)_{p=(a_2,...,a_{d^2+d})}$), which is invertible. Write
$$
M=
\begin{pmatrix}
A_p&B\\
C&D
\end{pmatrix},
$$

\noindent where $B$ is $\frac{d^2+5d}{2}\times 3$, $C$ is $\frac{d^2-3d+2}{2}\times \frac{d^2+5d}{2}$, and $D$ is $\frac{d^2-3d+2}{2}\times 3$.

Using the Guttman rank additivity formula and the Schur complement (see formula 0.9.2 of \cite{Schur}) we have that  $\operatorname{rank}(M)=\operatorname{rank}(A_p)+\operatorname{rank}(D-CA^{-1}B)$. Then 

$$
\operatorname{rank}M =\frac{d^2+5d}{2}
$$

\noindent it is equivalent to the vanishing of the entries of the Schur complement $D-CA^{-1}B$, which is of size $\frac{d^2-3d+2}{2} \times 3$. Than means $3\frac{d^2-3d+2}{2}$ determinantal equations.
\\

Now, we consider

$$
B=(b_1\ b_2\ b_3),
\qquad
D=(d_{ij}),
\qquad
C=\begin{pmatrix}
c_1\\
\vdots\\
c_{\frac{d^2-3d+2}{2}}
\end{pmatrix}
$$
\noindent where $b_j$ is  $\frac{d^2+5d}{2} \times 1$, $d_{ij} \in \C$; and $c_i$ is $1 \times \frac{d^2+5d}{2}$. Multiplying by the following matrix of determinant $1$ we obtain:
$$
\begin{pmatrix}
I_{\frac{d^2+5d}{2}}&0\\
-c_iA_p^{-1}&1
\end{pmatrix}
\begin{pmatrix}
A_p&b_j\\
c_i&d_{ij}
\end{pmatrix}=
\begin{pmatrix}
A_p&b_j\\
0&d_{ij}-c_iA_p^{-1}b_j
\end{pmatrix}.
$$

\noindent Since

$$
E_{ij}:
=
\det
\begin{pmatrix}
A_p&b_i\\
0&d_{ij}-c_iA_p^{-1}b_i
\end{pmatrix}
=
\det(A_p)(d_{ij}-c_iA_p^{-1}b_i).
$$

\noindent and $(D-CA_p^{-1}B)_{ij}=d_{ij}-c_iA_p^{-1}b_j$, we conclude that $\F_d^0$ is open in the algebraic variety:

$$\mathbb{V}\Big(E_{ij}:1 \leq i \leq \frac{d^2-3d+2}{2}, 1 \leq j \leq 3\Big).$$

For instance, when $d$ is odd, using the family obtained in Theorem 3 of \cite{alcantara2}, we consider the canonical monomial basis of $\mathbb{C}[y,z]_{\leq d+1}$, excluding $y^2$, $y^d z$, and $y^{d-1}$. By removing $\frac{(d-2)(d-1)}{2}$ elements from the basis of the codomain of $T(d+1,d^2+d+1)_p$, we obtain that the matrix $A_p$ described above is invertible. Therefore, we can determine the polynomial equations defining $\mathcal{F}_d^0$ on the open set $\det(A_p)\neq 0$.
\end{algorithm}

In what follows, we define another invariant that will be useful for our classification. In a subsequent work, we will explore the importance of this invariant in understanding the nature of the separatrices passing through the singularity for this type of foliations.

\begin{defin} Given $(a_2,...,a_{d^2+d}) \in \mathbb{A}^{d^2+d-1}$ such that $\Big\langle y-\sum_{i=2}^{d^2+d-1}a_i z^i,
z^{d^2+d}\Big\rangle$ is the local singular scheme of a foliation $\omega$. We will call \textbf{order of the singularity of $\omega$}, the natural number:

$$ord(\omega):=\min \{i: a_i \neq 0\}.$$
\end{defin}

Theorem 2 of \cite{alcantara2} establishes that this order is less than or equal to \(d\), the degree of the foliation. However, in the following theorem, we prove that the order must in fact be strictly less than \(d\). Among other things, this result will allows us to obtain a complete classification in degree 3, since in this case the order must be 2. This classification is presented in the next section.

\begin{theorem} \label{multiplicity} There are no foliations of degree $d\geq 3$ on $\mathbb{P}^2$
having a unique singularity of saddle-node or nilpotent type
with order $d$.
\end{theorem}

\begin{proof} According to  Proposition 4.2 of \cite{campillo2}, a homogeneous ideal $J \subset \C[x,y,z]$ of length $d^2+d+1$, is the singular scheme of a foliation on $\CP$ of degree $d$ if and only if there exist three linearly independent forms $A',B',C'$ in $H^0(\CP,J(d+1))$ without a common factor and three linearly independent linear forms  $L_{A'}, L_{B'}, L_{C'} \in H^0(\CP,\mathcal{O}(1))$ such that $J=\Big\langle A',B',C'  \Big\rangle$ and $L_{A'}A'+L_{B'}B'+L_{C'}C'=0$.  
\\

Let $\omega$ be a foliation on $\CP$ of degree $d$ with a unique singularity of multiplicity $1$ and order $d$, such that the singular scheme of the foliation has the following form:

$$I=\Big\langle y-z^d-\sum_{i=d+1}^{d^2+d}a_iz^i, z^{d^2+d+1} \Big\rangle \subset \C[y,z].$$ 

\noindent Up to a change of coordinates, we may normalize the coefficient of \(z^d\) to be equal to \(1\). We know that a polynomial $F(y,z) \in I$ if and only if $F(z^d+\sum_{i=d+1}^{d^2+d}a_iz^i,z) \in \langle z^{d^2+d+1} \rangle$ (see lemma \ref{lemma2}). 
 \\
 
 Let $\C[y,z]_k$ be  the vector space of homogeneous polynomials of degree $k$. In order to restrict the type of polynomials of degree \(d+1\) belonging to the ideal \(I\), we consider $W_d$, the vector space in $\C[y,z]_{\leq d+1}$ generate by the monomials: 

$$y,z^d,y\C[y,z]_k, \C[y,z]_{d+1}, \quad  \textrm{where $k \in \{1,...,d-1\}$};$$

\noindent and  $V$, the vector subspace of $\C[z]/\langle z^{d^2+d+1} \rangle$ generate by $z^d+\langle z^{d^2+d+1} \rangle, z^{d+1}+\langle z^{d^2+d+1} \rangle,...,z^{d^2+d}+\langle z^{d^2+d+1} \rangle$. We consider the vector space $W_d$ because a polynomial in $I$ can not have term $z^k$ for $1 \leq k < d$.
\\

 We have the following restriction the linear transformation 
$T(d+1,d^2+d+1)_{(a_{d+1},...,a_{d^2+d})}$:

\begin{align*}
T: W_d &\to V\\
f(y,z) &\mapsto f\Big(z^d+\sum_{i=d+1}^{d^2+d} a_iz^i,z\Big)+\langle z^{d^2+d+1} \rangle\\
\end{align*}

 \noindent We can see that the dimension of $W_d$ is  $\frac{d(d+3)}{2}+3$, and we use the canonical base of $W_d$: 
 
 $$\Big\{y,z^d,y^{j+1}z^{k-j},y^lz^{d+1-l}: 1\leq k \leq d-1, 0 \leq j \leq k, 0 \leq l \leq d+1\Big\}.$$

\noindent The lowest degree of $z$ in $T(f(y,z))$ will be call the order of the image of $f$ by $T$. If we consider the images by $T$ of the above basis, we obtain the following set of orders: 

$$\Big\{d,d(j+1)+k-j,dl+d+1-l: 1\leq k \leq d-1, 0 \leq j \leq k, 0 \leq l \leq d+1\Big\}.$$

\noindent This set has cardinality $\dim W_d-3=\frac{d(d+3)}{2}$, the only repeted orders are given by: 

\begin{align*}
\{T(y),T(z^d)\} &\quad \textrm{for order $d$},\\ 
\{T(yz),T(z^{d+1})\} &\quad \textrm{ for order $d+1$ and}\\ 
\{T(y^2),T(yz^d)\} &\quad \textrm{ for order $2d$}
\end{align*}

It is straightforward to see that elements of distinct orders in $V$ form a linearly independent set, and precisely we want the dimension of the image to be $dim W_d-3$. Because that means that the kernel has dimension $3$. Then the set: 

\begin{align*}
&\Big\{T(y),T(z^d),T(y^{j+1}z^{k-j}),T(y^lz^{d+1-l}): 1\leq k \leq d-1, 0 \leq j \leq k, 0 \leq l \leq d+1\Big\}\\
&-\Big\{T(z^d),T(z^{d+1}),T(yz^d)\Big\}.
\end{align*} 

\noindent is a $\C$-basis of the image of $T$. Then each element of $\{T(z^d),T(z^{d+1}),T(yz^d)\}$ admits a unique expression as a linear combination of the elements of this basis. This means that there exist unique $A_i, B_i,C_i$ homogeneous polynomial of degree $i$ in $\C[y,z]$ for $i=1,...,d+1$, such that:

\begin{align*}
T(z^d)&=-T(y)-\sum_{i=1}^{d-1}T\big(yA_i(y,z)\big)+T\big(A_{d+1}(y,z)\big)\\
T(z^{d+1})&=-T(yz)-\sum_{i=1}^{d-1}T\big(yB_i(y,z)\big)+T\big(B_{d+1}(y,z)\big)\\
T(yz^{d})&=-T(y^2)-\sum_{i=1}^{d-1}T\big(yC_i(y,z)\big)+T\big(C_{d+1}(y,z)\big).
\end{align*}

\noindent Therefore:

\begin{align*}
A&=y+z^d+\sum_{i=1}^{d-1}yA_i(y,z)+A_{d+1}(y,z),\\
B&=yz+z^{d+1}
+\sum_{i=1}^{d-1}yB_i(y,z)+yB_d(y,z),\\
C&=y^2+yz^d+\sum_{i=2}^{d-1}yC_i(y,z)+C_{d+1}(y,z),
\end{align*}

\noindent belong to the kernel of $T$ and we have: 

\begin{enumerate}
\item By construction, $A$ has neither a $z^{d+1}$ term nor a $yz^d$ term. Since the only ways to produce terms of degree $d+1$ and $2d$, respectively, are through the terms $yz$ and $y^2$, it follows that $A$ has neither a $yz$ term nor a $y^2$ term.
\item By construction, $B$ has neither a $z^{d}$ term nor a $yz^d$ term. Since the only ways to produce terms of degree $d$ and $2d$, respectively, are through the terms $y$ and $y^2$, it follows that $B$ has neither a $y$ term nor a $y^2$ term.
\item By construction, $C$ has neither a $z^{d}$ term nor a $z^{d+1}$ term. Since the only ways to produce terms of degree $d$ and $d+1$, respectively, are through the terms $y$ and $yz$, it follows that $C$ has neither a $y$ term nor a $yz$ term.

\end{enumerate}

\noindent Thus, we can write these polynomials as follows:

\begin{align*}
A&=y+z^d+\sum_{i=2}^{d-1}yA_i(y,z)+y^2D(y,z),\\
B&=yz+z^{d+1}
+\sum_{i=2}^{d-1}yB_i(y,z)+y^2E(y,z),\\
C&=y^2+yz^d+\sum_{i=2}^{d-1}yC_i(y,z)+y^2F(y,z),
\end{align*}

\noindent   where $D,E,F$ are homogeneous polynomials of degree $d-1$ in $\mathbb{C}[y,z]$. We conclude that $A,B,C$ are linearly independent in $W_d$, which means that they form a basis of the kernel. The homogenizations of these polynomials with respect to $x$ yield the following polynomials of degree $d+1$:

\begin{align*}
A(x,y,z)&=x^{d}y+xz^d+\sum_{i=2}^{d-1}x^{d-i}yA_i(y,z)+y^2D(y,z),\\
B(x,y,z)&=x^{d-1}yz+z^{d+1}
+\sum_{i=2}^{d-1}x^{d-i}yB_i(y,z)+y^2E(y,z),\\
C(x,y,z)&=x^{d-1}y^2+yz^d+\sum_{i=2}^{d-1}x^{d-i}yC_i(y,z)+y^2F(y,z),
\end{align*}

\noindent which belong to $H^0(\CP,I(d+1))$. Since $I$ is the singular scheme of a foliation, there exist: 

\begin{align*}
L_A&=\alpha_1 x+\alpha_2 y+\alpha_3 z,\\
L_B&=\beta_1 x+\beta_2 y+\beta_3z,\\
L_C&=\gamma_1 x+\gamma_2 y+\gamma_3 z
\end{align*}

\noindent linearly independet in $H^0(\CP,\mathcal{O}(1))$ such that $L_AA+L_BB+L_CC=0$. We can see that $\alpha_1=0,  (\beta_3=0)$, otherwise the term $x^{d+1}y$ in $L_AA$, ($z^{d+2}$ in $L_BB$) could not be eliminated.  
Suppose $\alpha_3=0$, then $y$ divides $L_B$, but this is a contradiction with linear independence. We can conclude that

\begin{align*}
L_A&=\alpha_2 y+z,\\
L_B&=-x+\beta_2 y,\\
L_C&=-\alpha_2 x+\gamma_2 y-\beta_2 z
\end{align*}

If $\alpha_2 \neq 0$, we have $C=\frac{L_A}{L_C}A+\frac{L_B}{L_C}B \in \langle A,B \rangle \mathcal{O}_{\C^2,0}$. Then  $\langle A,B \rangle \mathcal{O}_{\C^2,0}=I\mathcal{O}_{\C^2,0}$, therefore  the intersection index of $A$ and $B$ at $p=(0,0)$ is 
$$i_p(A,B)=d^2+d+1,$$
\noindent  and there are no other intersection points of $A$ and $B$ in $U_0$. Indeed, any such point would also be a zero of $C$ and hence a singular point for the foliation. Therefore, at infinity, that is, on $x=0$, $A$ and $B$ have $d$ intersection points. This implies that $A_{d+1}=y^2D(y,z)$ and 
$B_{d+1}=z^{d+1}+y^2E(y,z)$ have a common factor of degree $d$. However, this is impossible because $y$ does not divide $z^{d+1}+y^2E(y,z)$.
Therefore, $\alpha_2=0$. Thus, we have that there exist $\beta_2 \in \C$ and $\gamma_2 \in \C^*$ such that:

$$zA+(-x+\beta_2y)B+
(\gamma_2y-\beta_2z)C=0.$$

To finish the proof we prove the following:

\begin{prop}\label{proposition5}
Let $A,B,C\in\mathbb{C}[x,y,z]$ be homogeneous polynomials of degree
$d+1$, with $d\geq 3$, having no common factor. In the affine chart
$x=1$, suppose that they have the form
\[
A=
y+z^d+\sum_{i=2}^{d-1}yA_i(y,z)+y^2D(y,z),
\]
\[
B=
yz+z^{d+1}
+\sum_{i=2}^{d-1}yB_i(y,z)+y^2E(y,z),
\]
and
\[
C=
y\left(
y+z^d+\sum_{i=2}^{d-1}C_i(y,z)+yF(y,z)
\right),
\]
where $A_i,B_i,C_i$ are homogeneous polynomials of degree $i$ in $\C[y,z]$ and
$D, E, F$ are homogeneous of degree $d-1$ in $\C[y,z]$. Suppose that
\[
i_{(1:0:0)}(A,C)=d^2+d+1
\]
and that $(1:0:0)$ is the only intersection point of the curves
$A=0$ and $C=0$ in the affine chart $x=1$. Then there do not exist $\beta_2,\gamma_2\in\mathbb{C}$, with
$\gamma_2\neq0$, such that

$$zA+(-x+\beta_2y)B+
(\gamma_2y-\beta_2z)C=0.$$
\end{prop}

\begin{proof}
Since $B \in \langle A,C\rangle\mathcal{O}_{\C^2,0}$, we have $i_p(A,C)=d^2+d+1$, and there are no other intersection points of $A$ and $C$ in the affine chart $x=1$. Indeed, any such point would also be a zero of $B$ and hence a singular point of the foliation. Therefore, at infinity, that is, on $x=0$, $A$ and $C$ have $d$ intersection points. This implies that $A_{d+1}=y^2D$ and $C_{d+1}=y(yF+z^d)$ have a common factor of degree $d$. This common factor must be $yD$, since $y^2$ does not divide $C_{d+1}=y(yF+z^d)$. Hence, $yF+z^d=KD$, where $K$ is a linear form relatively prime to $y$.
\\

Assume that $\beta_2=0$. Then the relation

\begin{equation}\label{eq1}
zA+(-x+\beta_2y)B+(\gamma_2y-\beta_2z)C=0,
\end{equation}

\noindent becomes $zA-xB+\gamma_2yC=0$. Restricting \ref{eq1} to the line $x=0$, and using $A_{d+1}=y^2D$ and $C_{d+1}=yKD$, we obtain $zA_{d+1}+\gamma_2yC_{d+1}=0$. Hence, $zy^2D+\gamma_2y(yKD)=0$, and therefore $y^2D(z+\gamma_2K)=0$. Since $D\neq0$, it follows that $z+\gamma_2K=0$. As $\gamma_2\neq0$, we conclude that $K=-\frac{z}{\gamma_2}$.
\\
Since $i_p(A,C)=d^2+d+1$, we thus have that the intersection multiplicity
\begin{align*}
&i_p\big(y+z^d+\sum_{i=2}^{d-1}yA_i(y,z)+y^2D(y,z),y+z^d+\sum_{i=2}^{d-1}C_i(y,z)+yF(y,z)\big)=\\
&i_p\big(y+z^d+\sum_{i=2}^{d-1}yA_i(y,z)+y^2D(y,z),\sum_{i=2}^{d-1}yA_i(y,z)+y^2D(y,z)-\sum_{i=2}^{d-1}C_i(y,z)-yF(y,z)\big)=\\
&i_p\big(y+z^d+\sum_{i=2}^{d-1}yA_i(y,z)+y^2D(y,z),\sum_{i=2}^{d-1}yA_i(y,z)-\sum_{i=2}^{d-1}C_i(y,z)+y^2D+z^d+\frac{z}{\gamma_2}D\big)=\\
&d^2+1.
\end{align*}

Since $A=0$ is smooth at the origin, by the Implicit Function Theorem, the curve $A=0$ has a local branch at the origin given by $y=\varphi(z)$,
where $\varphi(z)=-z^d+\cdots$.
\\

\noindent Notice that, since we have the term $z^d$, the polynomial $D$ must contain the term $-\gamma_2z^{d-1}$ in order to eliminate the order $d$ when working with a branch of $A$. We therefore conclude that the second polynomial contains $y^2z^{d-1}$ as a term. Evaluating along the branch $\varphi(z)$, we obtain order $3d-1$. We now show that no other term of the polynomial
\[
\sum_{i=2}^{d-1}yA_i(y,z)-\sum_{i=2}^{d-1}C_i(y,z)+y^2D+z^d+\frac{z}{\gamma_2}D
\]
has this order.
\begin{enumerate}
\item The monomials of $yA_i$ are of the form $y^{j+1}z^{i-j}$, with $2\leq i\leq d-1$ and $0\leq j\leq i$. This monomial has order $d(j+1)+i-j$. Thus, if $d(j+1)+i-j=3d-1$, then $j(d-1)+i=2d-1$. If $j\leq1$, then $d\leq i$, a contradiction. If $j\geq2$, then $2d-1\geq2(d-1)+i$, which implies $i\leq1$, again a contradiction.
\item The monomials of $C_i$ are of the form $y^jz^{i-j}$, with $2\leq i\leq d-1$ and $0\leq j\leq i$. This monomial has order $jd+i-j$. Thus, if $jd+i-j=3d-1$, then $j(d-1)+i=3d-1$. If $j\leq2$, then $3d-1\leq2(d-1)+i$, which implies $d+1\leq i$, a contradiction. If $j\geq3$, then $jd+i-j=3d-1\geq3d+i-j$, which implies $j\geq i+1$, a contradiction.
\item The monomials of $D$ different from $z^{d-1}$ are of the form $y^jz^{d-1-j}$, with $1\leq j\leq d-1$. If $jd+d-1-j=3d-1$, then $j=2+\frac{2}{d-1}$. Since $d-1\geq2$, we have $2<j\leq3$. Hence $j=3$, and therefore $d=3$. But this is a contradiction, since $j\leq d-1$.
\end{enumerate}
\noindent We conclude that, since the order $3d-1$ cannot be eliminated, the intersection multiplicity satisfies
\[
d^2+1\leq3d-1.
\]
However, this is false for $d\geq3$. Therefore, $\beta_2\neq0$.
\\

Let $\beta_2\neq0$, since
\[
zA+(-x+\beta_2y)B+(\gamma_2y-\beta_2z)C=0,
\]
the intersection multiplicity at $p=(1:0:0)$ of the curves $A-\beta_2C$ and $B$ is
\begin{align*}
i_p\big(A-\beta_2C,B\big)&=i_p\big(A-\beta_2C,z(A-\beta_2C)-\gamma_2yC\big)\\
&=i_p\big(A-\beta_2C,yC\big)=(d+1)^2.
\end{align*}
\noindent This means that the curves $A-\beta_2C$ and $B$ have no intersection at infinity, that is, on the line $x=0$. On the other hand, on $x=0$ we have
\[
z(A_{d+1}-\beta_2C_{d+1})+\beta_2yB_{d+1}=0.
\]
Since $\beta_2\neq0$, it follows that $A_{d+1}-\beta_2C_{d+1}$ and $B_{d+1}$ have a common factor. This implies that there is an intersection point at infinity on $x=0$, which is a contradiction.

This concludes the proof of the proposition.
\end{proof}

We conclude that there do not exist polynomials $A$, $B$, and $C$ of degree $d+1$, defining a foliation of degree $d$ with a unique singularity whose linear part is nonzero and whose order is $d$.
\end{proof}

\begin{remark} For the case $d=2$, the intersection multiplicity of the curves $A$ and $C$ in proposition \ref{proposition5} is
\[
4d-1=d^2+d+1.
\]
This equality holds only in this case, and it is the unique degree for which such a foliation exists, with a unique saddle-node singularity of order $d$.
\end{remark}

\subsection{Foliations on $\CP$ of degree $3$ with a unique singular point}

In the case of foliations of degree $3$ with a unique singularity, the Milnor number is 13. In this setting, we can explicitly determine the polynomial equations defining the subvariety of foliations of this type for which the linear part at the singular point is nonzero. To this end, we will use again the theory of Schur complements of matrices.
\\

By Theorem \ref{multiplicity}, the curvilinear scheme, up to change of coordenates, of a foliation in $\mathcal{F}_3^0$ has the form: 
$$\left\langle y-z^2-\sum_{i=3}^{12}a_i z^i,
z^{13}\right\rangle,$$

\noindent then we will describe the subvariety of $\mathbb{A}^{10}$, which parametrizes foliations of degree $3$ with a unique singularity and nonzero linear part. We begin by proving some results that will subsequently allow us to provide the desired description of the subvariety.
\\

Let $p=(a_3,\dots,a_{12})\in \mathbb{A}^{10}$ and
\[g(z)=z^2+a_3z^3+\dots+a_{12}z^{12}=z^2u(z),\qquad u\in\C[z],\ u(0)=1 .\]

We have two elementary facts:

\begin{itemize}
\item  $ord(y^jz^k)=2j+k$, since $y^jz^k\mapsto z^{2j+k}u^j$ and $u(0)=1$.
\item $ord(fg)=ord(f)+ord(g)$.
\end{itemize}

Let \(A_p\) be the \(12\times 12\) block of the corresponding $13 \times 15$ matrix of $T(4,13)_p$; which rows are the coefficients of $z^0,\dots,z^{11}$ and columns are indexed by
\[
1,\ z,\ y,\ yz,\ yz^2,\ yz^3,\ y^2,\ y^2z,\ y^2z^2,\ y^3,\ y^3z,\ y^4,
\]
i.e.\ by $1$, $z$ and $y\,m$ with $m$ running over the ten monomials of degree $\le3$. %Put $\Delta(p):=\det A_p$.

\begin{lem}\label{lem:delta} With the above notation, we have the following
$$\det A_p=\det T(3,10)_p.$$ 

\noindent In particular $\det A_p=0$ if and only if there is $h\in\C[y,z]_{\le 3}$, $h\ne0$, with $ord(h)\ge10$.
\end{lem}

\begin{proof}
The columns $1$ and $z$ of $A_p$ are the unit vectors $e_0,e_1$. Every other column corresponds to some $y\,m$ and has image $g(z)\,m(g(z),z)$, divisible by $z^2$, hence has zeros in rows $z^0,z^1$. Thus $A_p=\begin{pmatrix}I_2&*\\0&A'\end{pmatrix}$, where $A'$ is the $10\times10$ block with rows $z^2,\dots,z^{11}$ and columns $y\,m$, therefore $\det A_p=\det A'$.

Let $m_1,\dots,m_{10}$ be the monomials of degree $\le3$, $g_i:=m_i(g(z),z)=\sum_{t\ge0}g_{i,t}z^t$, and $u=\sum_{s\ge0}u_sz^s$ with $u_0=1$. The entry of $A'$ in row $z^r$ ($2\le r\le11$) and column $y\,m_i$ is the coefficient of $z^r$ in $z^2u\,g_i$, i.e.\ the coefficient of $z^{r-2}$ in $u\,g_i$:
\[
A'_{r-2,\,i}=\sum_{t=0}^{r-2}u_{r-2-t}\,g_{i,t}\qquad(0\le r-2\le9).
\]
Let $T=(g_{i,t})_{t=0..9,\ i=1..10}$ be the matrix of $T(3,10)_p$ and $U=(u_{s-t})_{s,t=0}^{9}$ (with $u_k=0$ for $k<0$). Then $A'=U\,T$. The matrix $U$ is lower triangular with diagonal entries $u_0=1$, so $\det U=1$ and $\det A_p=\det A'=\det T$.
Finally $\det T=0$ iff $\ker T(3,10)_p\ne0$.
\end{proof}

\begin{theorem}\label{thm:main}
Let $p=(a_3,\dots,a_{12})\in \mathbb{A}^{10}$ be such that $T(3,12)_p$ is injective and $rank T(4,13)_p\le12$. Then $\det A_p \neq0$.
\end{theorem}

\begin{proof}
Fix such a $p$ and suppose, for contradiction, that $\det A_p=0$. We will construct $13$ polynomials $f_0,\dots,f_{12}\in\C[y,z]_{\le4}$ with $ord(f_k)=k$ for $k=0,\dots,12$. This forces $rank T(4,13)_p=13$.

To begin with, we can observe that the monomials
\[
f_0=1,\ f_1=z,\ f_2=z^2,\ f_3=z^3,\ f_4=z^4,\ f_5=yz^3,\ f_6=y^2z^2,\ f_7=y^3z,\ f_8=y^4
\]
lie in $\C[y,z]_{\le4}$ and satisfy $ord(f_k)=k$. 
\\

Let $M$ be the $12\times10$ matrix of $T(3,12)_p$ (its rows are the coefficients of $z^0,\dots,z^{11}$ and its columns correspond to the ten monomials of degree $\le3$). Since $T(3,12)_p$ is injective, $rank M=10$. Perform Gaussian elimination on the columns of $M$. Then, each of the ten columns has a first nonzero entry, and these first nonzero entries lie in the ten pairwise distinct rows $n_1<n_2< \dots<n_{10}\le 11$.
\\

Let $h_j\in\C[y,z]_{\le3}$ be the corresponding polynomial of $n_j$, for $j=1,...,10$. The $i$-th column is the coefficient vector of $h_i(g(z),z)\bmod z^{12}$, whose first nonzero entry is in row $n_i\le11$; hence $ord(h_i)=n_i$. 
\\

In fact, we can note that every nonzero $h\in\C[y,z]_{\le3}$ satisfies $ord(h) \in\{n_1,\dots,n_{10}\}$. Indeed, write $h=\sum c_ih_i$ and let $i_0$ be the smallest index with $c_{i_0}\ne0$; since the $n_i$ are strictly increasing, then $ord(h)=n_{i_0}$. The monomials $1,z,z^2,z^3,yz^2,y^2z,y^3\in\C[y,z]_{\le3}$ have orders $0,1,2,3,4,5,6$. Therefore, these seven values must be among the $n_i$. Since the $n_i$ are increasing then $n_1,\dots,n_7 \in \{0,\dots,6\}$ and $7\le n_8<n_9<n_{10}\le 11$.
\\

From $\det A_p =0$, we have that there exists $h_0\in \C[y,z]_{\le3} \setminus\{0\}$ with $ord(h_0)\ge10$. Since $ord(h_0)$ is some $n_i$, which is greater or equal than $10$, it is one of $n_9,n_{10}$. This forces $n_{10}\in\{10,11\}$.

We have nine polynomials $f_k$, for $k=0,...,8$. We will now construct the remaining four to prove that the rank of \(T(4,13)_p\) is 13.
We use only $h_8,h_9,h_{10}$ and the fact that,
\[
ord(h_i)=n_i,\qquad ord(zh_i)=n_i+1,\qquad ord(yh_i)=n_i+2,
\]
where $h_i,\ zh_i,\ yh_i \in\C[y,z]_{\le4}$ because $\deg h_i\le3$. Based on the above, we have only the following possibilities.

\emph{Case $n_{10}=10$.} Then $7\le n_8<n_9\le9$, so $n_8\in\{7,8\}$. Take
\[
f_{10}=h_{10},\qquad f_{11}=zh_{10},\qquad f_{12}=yh_{10},
\]
of orders $10,11,12$, and
\[
f_9=\begin{cases} yh_8 & \text{if } n_8=7,\\ zh_8 & \text{if } n_8=8,\end{cases}
\]
which has order $9$ in both cases.

\emph{Case $n_{10}=11$.} Take $f_{11}=h_{10}$ and $f_{12}=zh_{10}$, of orders $11,12$. Here $7\le n_8<n_9\le10$, so $(n_8,n_9)$ is one of six pairs, and we choose $f_9,f_{10}$ as follows.
\begin{center}
\begin{tabular}{c|c|c|l}
$(n_8,n_9)$ & $f_9$ (order $9$) & $f_{10}$ (order $10$) \\\hline
$(7,8)$  & $yh_8$ & $yh_9$ \\
$(7,9)$  & $h_9$  & $zh_9$ \\
$(7,10)$ & $yh_8$ & $h_9$ \\
$(8,9)$  & $zh_8$ & $yh_8$ \\
$(8,10)$ & $zh_8$ & $h_9$  \\
$(9,10)$ & $h_8$  & $h_9$  
\end{tabular}
\end{center}

In all cases we have obtained $f_9,f_{10},f_{11},f_{12}\in\C[y,z]_{\le4}$ with $ord f_k=k$. Therefore $rank T(4,13)_p=13$.

\end{proof}

Given $p=(a_3,...,a_{12}) \in \mathbb{A}^{10}$. For $i=1,2,3$, let \(E_i\) be the determinat of the \(13\times 13\) block of  the matrix of $T(4,13)_p$; which rows are the coefficients of $z^0,\dots,z^{12}$ and columns are indexed by
\[
1,\ z,\ y,\ yz,\ yz^2,\ yz^3,\ y^2,\ y^2z,\ y^2z^2,\ y^3,\ y^3z,\ y^4, z^{i+1}.
\]

\begin{theorem} \label{degree3} For foliations of degree $3$ on $\CP$ with a unique singularity with multiplicity $1$, we have the following description of the algebraic variety:

\begin{align*}
\F_3^0&=\Big\{(a_3,...,a_{12}) \in \mathbb{A}^{10} : \textrm{$\langle y-z^2-\sum_{i=3}^{12}a_i z^i,
z^{13}\rangle$ is the singular scheme of a foliation}  \Big\}\\
& =\mathbb{V}(E_1,E_2,E_3) \cap \{p \in \mathbb{A}^{10}: \det(T(3,10)_p)\neq 0\},
\end{align*}

\noindent and  $\dim \F_3^0= 7$.
\end{theorem}

\begin{proof} We begin by prooving the following:
$p=(a_3,...,a_{12}) \in \mathbb{A}^{10}$ defines the singular scheme of a foliation of degree $3$ in and only if $nul(T(3,12)_p)=0$ and $rank(T(4,13)_p) \leq 12$.  If $p=(a_3,...,a_{12}) \in \mathbb{A}^{10}$ defines the singular scheme of a foliation of degree $3$, by Lemma 1, we have $nul(T(3,12)_p)=0$ and $rank(T(4,13)_p) \leq 12$. On the other hand, suppose that $nul(T(3,12)_p)=0$ and $rank(T(4,13)_p) \leq 12$. Note that the vectors $T(1,6)_p(1)=1+\langle z^6 \rangle$, $T(1,6)_p(z)=z+\langle z^6 \rangle$ and $T(1,6)_p(y)=z^2+a_3z^3+a_4z^4+a_5z^5+\langle z^6 \rangle$ are linearly independent in the image, then $nul(T(1,6)_p)=0$. By Corollary 4.9 b) of \cite{campillo2} we have that  $p=(a_3,...,a_{12}) \in \mathbb{A}^{10}$ defines the singular scheme of a foliation.
\\

Let $p=(a_3,...,a_{12}) \in \mathbb{A}^{10}$ and, as above, let \(A_p\) be the \(12\times 12\) block of the corresponding $13 \times 15$ matrix $M$ of $T(4,13)_p$; which rows are the coefficients of $z^0,\dots,z^{11}$ and columns are indexed by $1,\ z,\ y,\ yz,\ yz^2,\ yz^3,\ y^2$, $\ y^2z,\ y^2z^2,\ y^3,\ y^3z,\ y^4$.

Write
$$
M=
\begin{pmatrix}
A_p&B\\
C&D
\end{pmatrix},
$$

\noindent where $B$ is $12\times3$, $C$ is $1\times12$, and $D$ is $1\times3$.

If $p \in \F_3^0$, then $A_p$ is invertible. Using the Guttman rank additivity formula and the Schur complement (see formula 0.9.2 of \cite{Schur}) we have that  $\operatorname{rank}(M)=\operatorname{rank}(A_p)+\operatorname{rank}(D-CA^{-1}B)$. Then 

$$
\operatorname{rank}M\leq12
$$

\noindent it is equivalent to the vanishing of the three entries of the Schur complement $D-CA^{-1}B$ of size $1 \times 3$.
\\

Now, we consider

$$
B=(b_1\ b_2\ b_3),
\qquad
D=(d_1\ d_2\ d_3),
$$
\noindent with this, we recover the polynomials \(E_1, E_2, E_3\) that we defined previously.:
$$
E_i=
\det
\begin{pmatrix}
A_p&b_i\\
C&d_i
\end{pmatrix},
$$
\noindent for $i=1,2,3$. Multiplying by the matrix of determinant $1$,
$$
\begin{pmatrix}
I_{12}&0\\
-CA_p^{-1}&1
\end{pmatrix}
\begin{pmatrix}
A_p&b_i\\
C&d_i
\end{pmatrix},
$$
we obtain
$$
E_i
=
\det
\begin{pmatrix}
A_p&b_i\\
0&d_i-CA_p^{-1}b_i
\end{pmatrix}
=
\det(A_p)(d_i-CA_p^{-1}b_i).
$$
Since $D-CA_p^{-1}B=\big( d_1-CA_p^{-1}b_1 \quad d_2-CA_p^{-1}b_2  \quad d_3-CA_p^{-1}b_3\big)$. With this we have that $\F_3^0 \subset \mathbb{V}(E_1,E_2,E_3) \cap \{p \in \mathbb{A}^{10}: \det(T(3,10)_p)\neq 0\}$. 
\\

The reverse inclusion is clear: if $p \in  \mathbb{V}(E_1,E_2,E_3) \cap \{p \in \mathbb{A}^{10}: \det(T(3,10)_p)\neq 0\}$ then the matrix $T(4,13)_p$ has rank $12$. And since $Ker(T(3,12)_p) \neq \{0\}$ implies  $Ker(T(3,10)_p) \neq \{0\}$, because a cubic passing through 12 points also passes through 10 points, we have the result.

\end{proof}

Now we present a particular example that allows us to explicitly construct families of foliations of degree 3 with a unique singularity whose linear part is nonzero. This family has dimension three and generalizes the one given in page 197 of \cite{alcantaraPantaleon}. These are foliations whose singularity is a node-saddle. Examples of foliations with nilpotent singularities can be found in Theorem 3 of \cite{alcantara2}.

\begin{example} \label{familydegree3}  The curvilinear scheme:

$$I=\left\langle y-z^2-a_5z^5-a_8z^8-a_{11}z^{11}, z^{13} \right\rangle \subset \C[y,z] $$

\noindent is the singular scheme of a foliation of degree $3$ on $\CP$ with a unique singularity if and only if: $3a_5^4-2a_5^2a_8-a_8^2+a_5a_{11}=0$ and $a_5(2a_5^2-a_8)(3a_5^2-a_8) \neq 0$.
\end{example}

\begin{proof} 

For $p=(0,0,a_5,0,0,a_8,0,0,a_{11},0)\in\mathbb A^{10}$, we will prove that \(E_1(p)=E_3(p)=0\). Let $g(z)=z^2+a_5z^5+a_8z^8+a_{11}z^{11}$. Since every exponent appearing in \(g(z)\) is congruent to \(2\) modulo \(3\), every exponent appearing in \(z^kg(z)^r\) is congruent to \(k+2r\) modulo \(3\).

Expanding the determinant \(E_i\) first along the column corresponding to \(z^{i+1}\), and then along the columns corresponding to \(1\) and \(z\), shows that, up to sign, \(E_i\) is the determinant of the \(10\times10\) matrix whose columns correspond to
\[
g,\ zg,\ z^2g,\ z^3g,\ g^2,\ zg^2,\ z^2g^2,\ g^3,\ zg^3,\ g^4,
\]
and whose rows correspond to the coefficients of \(z^2,\ldots,z^{12}\), with the row corresponding to \(z^{i+1}\) deleted.

For \(E_1\), after deleting the row corresponding to \(z^2\), only three rows remain with exponents congruent to \(2\) modulo \(3\), namely those corresponding to \(z^5,z^8,z^{11}\). However, the four columns corresponding to $g,\quad z^3g,\quad zg^2,\quad g^4$,
are supported only on these three rows. Hence these columns are linearly dependent, and \(E_1(p)=0\).

Similarly, for \(E_3\), after deleting the row corresponding to \(z^4\), only two rows remain with exponents congruent to \(1\) modulo \(3\), namely those corresponding to \(z^7,z^{10}\). The three columns corresponding to $z^2g,\quad g^2,\quad zg^3$
are supported only on these two rows. Thus, they are linearly dependent, and \(E_3(p)=0\). By a direct computation, we obtain:

\begin{align*}
E_2(p)=&(-2a_5^2+a_8)(3a_5^2-a_8) \left(3a_5^4-2a_5^2a_8+a_5a_{11}-a_8^2\right)\\
\det(T(3,10)_p)=&a_5(2a_5^2-a_8)(3a_5^2-a_8).
\end{align*}

\noindent Therefore $I=\langle y-z^2-a_5z^5-a_8z^8-a_{11}z^{11}, z^{13} \rangle$ is the singular scheme of a foliation of degree $3$ if and only if $3a_5^4-2a_5^2a_8-a_8^2+a_5a_{11}=0$ and $a_5(2a_5^2-a_8)(3a_5^2-a_8) \neq 0$.
\\

Now, given $a_3,a_8,a_{11} \in \C$, such that that $3a_5^4-2a_5^2a_8-a_8^2+a_5a_{11}=0$ and $a_5(2a_5^2-a_8)(3a_5^2-a_8) \neq 0$ we will construct the corresponding foliation. Consider the polynomials:

\begin{align*}
F(y,z)=&
y\big((a_5^4-a_5^2a_8+a_5a_{11}-a_8^2)y^3+(7a_5^3-6a_5a_8+a_{11})z^3\big)\\
&+(-4a_5^3+5a_5a_8-a_{11})zy^2+(3a_5^2-a_8)z^2+(-3a_5^2+a_8)y,\\[6pt]
G(y,z)=&
z \big(-a_5^3y^3+(2a_5^2-a_8)z^3 \big)\\
&+(-5a_5^2+2a_8)z^2y+(3a_5^2-a_8)y^2,
\end{align*}

\noindent we can see that

\begin{align*}
&F(a_2z^2+a_5z^5+a_8z^8+a_{11}z^{11},z)=z^{14}F_1(z)\\
&G(a_2z^2+a_5z^5+a_8z^8+a_{11}z^{11},z)=z^{13}G_1(z),
\end{align*}

\noindent therefore $F,G \in I$. Let $\alpha_1=a_5^4-a_5^2a_8+a_5a_{11}-a_8^2$, $\alpha_2=7a_5^3-6a_5a_8+a_{11}$, $\beta_1=-a_5^3$ and 
$\beta_2=2a_5^2-a_8$, then $\alpha_1\beta_2-\beta_1\alpha_2=
(3a_5^2-a_8)
\left(
3a_5^4+a_5a_{11}-2a_5^2a_8-a_8^2
\right)=0$. With $F$ and $G$ we construct the following three homogeneous polynomials of degree four:

\[
\begin{aligned}
A=&
\Big[(2a_5^2-a_8)(-4a_5^3+5a_5a_8-a_{11}) -(7a_5^3-6a_5a_8+a_{11})(-5a_5^2+2a_8)\Big]y^2z^2\\
&+(3a_5^2-a_8)x\Big[(2a_5^2-a_8)z^3 -(7a_5^3-6a_5a_8+a_{11})y^3\Big]-(2a_5^2-a_8)(3a_5^2-a_8)x^2yz, \\
B=&(7a_5^3-6a_5a_8+a_{11})\Big[
(2a_5^2-a_8)z^4-a_5^3zy^3+(-5a_5^2+2a_8)xz^2y
+(3a_5^2-a_8)x^2y^2
\Big],\\
C=&-(2a_5^2-a_8)\Big[
(7a_5^3-6a_5a_8+a_{11})z^3y +(a_5^4-a_5^2a_8+a_5a_{11}-a_8^2)y^4\\
&+(-4a_5^3+5a_5a_8-a_{11})xzy^2+(3a_5^2-a_8)x^2z^2
-(3a_5^2-a_8)x^3y
\Big].
\end{aligned}
\]

\noindent These polynomials satisfy:

\[
xA(x,y,z)+yB(x,y,z)+zC(x,y,z)=0,
\]
then the homogeneous 1-form
\[
\omega=A(x,y,z)\,dx+B(x,y,z)\,dy+C(x,y,z)\,dz
\]
defines a holomorphic foliation of degree three on \(\mathbb P^2\), with a unique singular point, which is of node-saddle type.

\end{proof}

\subsection{Foliations on $\CP$ of degree $4$ with a unique singular point}

In this final section, we present an example of a foliation on $\CP$ of degree $4$, with a unique singularity with non-zero linear part. This example is particularly interesting because, until now, general examples were known only in odd degrees, whereas no examples were available in even degrees other than degree two.

The construction of this foliation was inspired by an example obtained through an interaction with Claude (Anthropic) over a field of characteristic seven. By using tools such as Hensel's lemma, we were then able to obtain the corresponding example in characteristic zero. Unfortunately, the resulting ideal does not seem to exhibit any discernible pattern that could be exploited to construct examples in even degrees greater than or equal to six.

\begin{prop} \label{exampledegree4}
The ideal $I=\langle y-g(z),z^{21} \rangle$, where
\[
\begin{aligned}
g(z)&=z^3+z^7+z^8+z^9+z^{10}+z^{11}+z^{12}\\
&-\big(+8z^{13}+31z^{14}+104z^{15}+276z^{16}+733z^{17}+1834z^{18}+4508z^{19}+10723z^{20} \big),
\end{aligned}
\]
is the local scheme of a foliation on $\CP$ of degree $4$ having a unique singular point, which is of saddle-node type. The foliation is defined by the homogeneous \(1\)-form
\[
\omega=A^h(x,y,z)\,dx+B^h(x,y,z)\,dy+C^h(x,y,z)\,dz,
\]
where \(A^h,B^h,C^h\) are the homogenizations of degree \(5\)
with respect to the variable \(x\) of the following polynomials:

\begin{align*}
A =& -88y^5 - 57y^4z + 71y^3z^2 - 142y^2z^3 + 93yz^4 - 3z^5 + 73y^4 - 22y^3z - y^2z^2 \\
&- 18yz^3 - z^4 + 138y^3 - 93y^2z + 3yz^2 + 18y^2 + yz, \\
B =& -243y^4z - 164y^3z^2 + 194y^2z^3 - 387yz^4 + 247z^5 + 88y^4 + 259y^3z - 126y^2z^2 \\
&+ 138yz^3 - 144z^4 - 73y^3 + 405y^2z - 247yz^2 + 18z^3 - 138y^2 + 144yz - 18y, \\
C =& 243y^5 + 164y^4z - 194y^3z^2 + 387y^2z^3 - 247yz^4 - 202y^4 + 55y^3z + 4y^2z^2 \\
&+ 51yz^3 + 3z^4 - 383y^3 + 248y^2z + z^3 - 51y^2 - 3yz - y.
\end{align*}
\end{prop}

\begin{proof} A direct computation shows that $A+yB+zC=0$. Homogenizing this identity yields the Euler relation $xA^h+yB^h+zC^h=0$. Direct substitution gives

\begin{align*}
A(g(z),z)
&=25067z^{22}+O(z^{23})\\
B(g(z),z)
&=-450365z^{21}+O(z^{22}),\\
C(g(z),z)
&=-25067z^{21}+O(z^{22}).
\end{align*}

\noindent Therefore $A,B,C \in I$ and using the free software Singular with the lexicographic monomial ordering, we obtain the following reduced Gröbner basis of the ideal generated by
 \(B\) and \(C\):

\begin{verbatim}
> ring r=0,(y,z),lp;
> poly B = -243*y4z - 164*y3z2 + 194*y2z3 - 387*yz4 + 247*z5 + 88*y4 +
 259*y3z - 126*y2z2 +138*yz3 - 144*z4 - 73*y3 + 405*y2z - 
 247*yz2 + 18*z3 - 138*y2 + 144*yz - 18*y;
> poly C = 243*y5 + 164*y4z - 194*y3z2 + 387*y2z3 - 
247*yz4 - 202*y4 + 55*y3z +4*y2z2 + 51*yz3 + 3*z4 - 383*y3 + 248*y2z + 
z3 - 51*y2 - 3*yz - y;
> ideal I=(B,C);
> ideal Gr = groebner(I);
> Gr;
Gr[1]=z21
Gr[2]=y+10723z20+4508z19+1834z18+733z17+276z16+104z15+31z14+8z13
-z12-z11-z10-z9-z8-z7-z3
\end{verbatim}

\noindent Thus, the Gröbner basis is precisely the ideal \(I\). Following Theorem 1 of \cite{alcantara2} we conclude that the foliation $\omega$ on $\CP$ of degree $4$, given locally by $\omega_0=Bdy+Cdz$ has a unique singular point at $(1:0:0)$. The linear part of the foliation is $-18ydy-ydz$, therefore the singularity is a saddle-node. 

\end{proof}

Consider the above point 

$$p = (0, 0, 0, 1, 1, 1, 1, 1, 1, -8, -31, -104, -276, -733, -1834, -4508, -10723) \in \mathbb{A}^{17}$$ 

\noindent  and take the $18\times 18$ submatrix of $T(5,21)_p$ obtained from the canonical bases of $\mathbb{C}[y,z]_{\le 5}$ (the $21$ monomials $y^i z^j$ with $i+j\le 5$) and the canonical basis of $\mathbb{C}[z]/ \langle z^{21} \rangle$ (the $21$ elements $1, z, \dots, z^{20}$), where we have excluded:

\begin{enumerate}
  \item from the basis of $\mathbb{C}[y,z]_{\le 5}$, the monomials $y$, $y^2$, $y^2 z$;
  \item from the basis of $\frac{\mathbb{C}[z]}{\langle z^{21} \rangle}$, the elements $z^{17}$, $z^{19}$, $z^{20}$.
\end{enumerate}

\noindent We can see that the determinant of this matrix at $p$ is not zero. Then in the open set given by que complement of the determinant of this matrix we can apply the complement the Schur to have de variety $\mathcal{F}_4^0$, which is defined by 9 polynomials (see algorithm \ref{algorithm}).

\section*{Acknowledgements}

The author would like to thank Dominique Cerveau for fruitful
discussions on the subject of this paper and Manuel Cruz for pointing out an appropriate reference for the version of Hensel's lemma used in this work.

\end{document}